%% file: main.tex
\documentclass[letterpaper, 11pt]{article}
\usepackage[utf8]{inputenc}

\usepackage{geometry}
\usepackage{amssymb}
\usepackage{amsmath}
\usepackage{breqn}
\usepackage{enumerate}
\usepackage{algorithm}
\usepackage{algpseudocode}
\usepackage{amsthm}
\usepackage{hyperref}
\usepackage{tikz}
\usetikzlibrary{shapes.geometric, arrows, backgrounds, scopes, positioning}
\usepackage{graphicx}
\usepackage{subcaption}
\usepackage{colortbl}
\usepackage{cite}
\usepackage{authblk}
\usepackage{multirow}

\providecommand{\norm}[1]{\ensuremath{\left\lVert#1\right\rVert }}

\def\R{\mathbb{R}}
\def\A{A^TA}

\def\E{\mathbb{E}}
\def\Ezt{\mathbb{E}_{\zeta_t}}
\providecommand{\expect}[1]{\E_t\ensuremath{\left[#1\right]}}
\providecommand{\fixedexp}[1]{\E\ensuremath{\left[#1\right]}}
\providecommand{\serverexp}[1]{\Ezt\ensuremath{\left[#1\right]}}
\newtheorem{theorem}{\bfseries Theorem}
\newtheorem{lemma}{\bfseries Lemma}
\newtheorem{assumption}{\bfseries Assumption}

\title{
Robustness of Iteratively Pre-Conditioned Gradient-Descent Method: The Case of Distributed Linear Regression Problem}

\author{Kushal Chakrabarti$^\star$, Nirupam Gupta$^\dagger$, and Nikhil Chopra$^\star$
\thanks{$^\star$ University of Maryland, College Park, Maryland 20742, U.S.A. \\
$^\dagger$ Georgetown University, Washington, DC 20057, U.S.A. \\
Emails: {\em kchak@terpmail.umd.edu}, {\em nirupam115@gmail.com}, and {\em nchopra@umd.edu}}%
}

\date{}

\begin{document}

\maketitle

\input{abstract}
\input{intro}
\input{main_results}
\input{proofs}

\input{simulations}
\input{summary}

\newpage
\bibliographystyle{unsrt}
\bibliography{refs}


\addtolength{\textheight}{-12cm}

\end{document}

%% file: abstract.tex
\begin{abstract}

This paper considers the problem of multi-agent distributed linear regression in the presence of system noises.
In this problem, the system comprises multiple agents wherein each agent locally observes a set of data points, and the agents' goal is to compute a linear model that best fits the collective data points observed by all the agents. We consider a server-based distributed architecture where the agents interact with a common server to solve the problem; however, the server cannot access the agents' data points. 
We consider a practical scenario wherein the system either has {\em observation noise}, i.e., the data points observed by the agents are corrupted, or has {\em process noise}, i.e., the computations performed by the server and the agents are corrupted.
In noise-free systems, the recently proposed distributed linear regression algorithm,
named the Iteratively Pre-conditioned Gradient-descent (IPG) method, has been claimed to converge faster than related methods. In this paper, we study the {\em robustness} of the IPG method, against both the {\em observation noise} and the {\em process noise}. We empirically show that the robustness of the IPG method compares favorably to the state-of-the-art algorithms. 

\end{abstract}

%% file: intro.tex
\section{INTRODUCTION}
\label{sec:intro}

This paper considers the problem of multi-agent distributed linear regression in the presence of additive system noises, namely the {\em observation noise} and {\em process noise}. The nomenclature {\em distributed} refers to the data being distributed amongst multiple agents~\cite{azizan2019distributed, bertsekas1989parallel}. In this problem, the goal is to design an algorithm that -
\begin{itemize}
    \item allows the agents to compute an optimal linear regression model for their collective data (see Equation~\eqref{eqn:opt_1} below), 
    \item does not require any individual agent to share its local data points.
\end{itemize}
Distributed regression problems have received significant attention in recent years, see~\cite{li2014scaling, yang2019federated}, owing to the increased availability of data, handheld devices, and ubiquitous communication networks.\\

\tikzstyle{master} = [rectangle, rounded corners, minimum width=1.7cm, minimum height=1cm,text centered, text width=1cm, draw=black, fill=blue!30]
\tikzstyle{machine} = [rectangle, minimum width=1cm, minimum height=1cm,text centered, text width=2cm, draw=black, fill=blue!10]
\tikzstyle{dots} = [circle, inner sep=0pt,minimum size=2pt, draw=black, fill=blue!50!cyan]
\tikzstyle{arrow} = [thick,<->,>=stealth]

\begin{figure}[thpb]  
\centering
\begin{tikzpicture}[node distance = 1.5cm, auto]
    \node (master) [master] {Server};
    \node (m/c2) [machine, below  = of master] {Agent 2\\ $(A^2,B^2)$};
    \node (m/c1) [machine, left = of m/c2, xshift=1cm] {Agent 1\\ $(A^1,B^1)$};
    \node (d1) [dots, right = of m/c2, xshift=-0.8cm] {};
    \node (d2) [dots, right = of d1, xshift=-1.4cm] {};
    \node (d3) [dots, right = of d2, xshift=-1.4cm] {};
    \node (m/c3) [machine, right = of d3, xshift=-0.8cm] {Agent m\\ $(A^m,B^m)$};
    
    \draw[arrow] (master) -- (m/c1);
    \draw[arrow] (master) -- (m/c2);
    \draw[arrow] (master) -- (m/c3);
\end{tikzpicture}
\caption{System architecture.}
\label{fig:sys}
\end{figure}
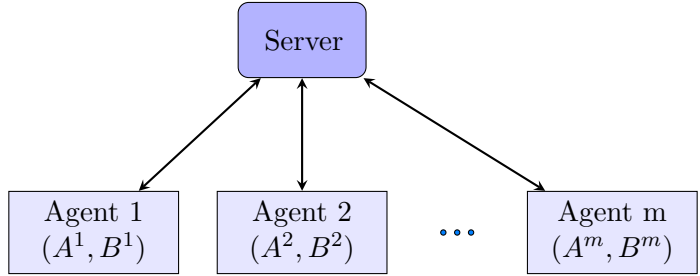


Specifically, we consider a {\em server-based} system architecture comprising $m$ agents and one server, shown in Fig.~\ref{fig:sys}. The agents can only interact with the server, and the overall system is assumed synchronous. Each agent $i\in \{1, \ldots, \, m\}$ has $n_i$ {\em local} data points, represented by an {\em input matrix} $A^i \in \R^{n_i \times d}$ and an {\em output vector} $b^i \in \R^{n_i}$.
The agents \underline{do not} share their individual local data points with the server. For each agent $i$, we define a {\em local cost function} $F^i: \R^d \to \R$ such that for a given {\em regression parameter} $x \in \R^d$,
\begin{align}
     F^i(x) = \frac{1}{2} \norm{A^i x - b^i}^2, \label{eqn:dist_opt}
\end{align}
where $\norm{\cdot}$ denotes the Euclidean norm.
A distributed regression algorithm prescribes instructions for coordination of the server and the agents to compute an optimal regression parameter $x^* \in \R^d$ such that
\begin{align}
  x^* \in X^* = \arg \min_{x \in \R^d} \, \sum_{i = 1}^m F^i(x). \label{eqn:opt_1}
\end{align}
In this paper, we refer to the above server-based system architecture as {\em server-agent network}. Common applications of the above problem include distributed state estimation, hypothesis testing~\cite{zhu2018linear}, supervised machine learning problems such as the supply chain demand forecasting~\cite{carbonneau2008application}, predicting online user input actions~\cite{canny2013method}, sparse linear solver selection~\cite{bhowmick2006application}.\\

There are numerous distributed algorithms available for solving the above regression problem, for example, see~\cite{chakrabarti2020iterative,azizan2019distributed} and references therein. A few notable ones are the gradient-descent (GD) method~\cite{bertsekas1989parallel}, Nesterov's accelerated gradient-descent method (NAG)~\cite{nesterov27method}, the heavy-ball method (HBM)~\cite{polyak1964some}, the accelerated projection-consensus method (APC)~\cite{azizan2019distributed}, the quasi-Newton BFGS method~\cite{kelley1999iterative}, and a recently proposed iteratively pre-conditioned gradient-descent (IPG) method~\cite{chak2020ipc}. These distributed algorithms are iterative wherein the server maintains an estimate of a solution defined by~\eqref{eqn:opt_1}, which is updated iteratively using the gradients of the individual agents' local cost functions defined in~\eqref{eqn:dist_opt}. In an ideal scenario with no noise, these algorithms converge to an optimal regression parameter defined in~\eqref{eqn:opt_1}.\\

Practical systems, however, inevitably suffers from uncertainties or {\em noise}~\cite{gold1966effects, holi1993finite}. Specifically, we consider two types of additive system noises, 1) {\em observation noise}, and 2) {\em process noise}. The {\em observation noise}, as the name suggests, models the uncertainties in the local data points observed by the agents~\cite{zhu2004class}. The {\em process noise} models the uncertainties or jitters in the computation process due to hardware failures, quantization errors, or noisy communication links~\cite{holi1993finite}. In this paper, we empirically show (for the first time) that the approximation error, or the {\em robustness}, of the IPG method in the presence of the above system noises \underline{compares favorably} to all the other aforementioned prominent distributed algorithms.\\

We note, however, that the prior work on the formal convergence of 
the IPG method only considers a Utopian setting wherein the system is free from noise~\cite{chak2020ipc}. Therefore, besides empirical results, we also present formal analyses on the IPG method's robustness when negatively impacted by additive system noises. Robustness analyses of some of the other aforementioned distributed algorithms are in~\cite{devolder2014first,sebbouh2020convergence}.

\subsection{Summary of Our Contributions}
\label{sub:contri}
The IPG method was first proposed in~\cite{chak2020ipc} and claimed to converge at a faster rate than the aforementioned state-of-the-art methods in the noise-free case~\cite{chakrabarti2020iterative}. Here, we make two key contributions, summarized as follows.

\begin{itemize}
\setlength\itemsep{0.5em}
    \item {\bf Theory :} In Section~\ref{sec:results}, we theoretically characterize the robustness guarantees of the IPG method against both the observation and the process noises.
    
    \item {\bf Experiments:} In Section~\ref{sec:exp}, we empirically show the improved robustness of the IPG method, in comparison to the state-of-the-art algorithms.
\end{itemize}

%% file: main_results.tex
\section{ROBUSTNESS OF THE IPG METHOD}
\label{sec:results}

\begin{algorithm}
  \caption{The IPG method for the noise-free case~\cite{chak2020ipc}.}\label{algo_1}
  \begin{algorithmic}[1]
    \State The server initializes $x(0) \in \R^d$, $K(0) \in \R^{d \times d}$, $\alpha > 0$, $\delta > 0$.
    \For{\text{each iteration $t=0, \, 1, \, 2, \ldots $}}
      \State The server broadcasts $x(t)$ and $K(t)$ to all the agents.
      \State Each agent $i  \in \{1, \ldots, \, m\} $ sends to the server a gradient $g^i(t)$, defined as
      \begin{align}
        g^i(t) = \nabla F^i(x(t)) = \left( A^i \right)^T \, \left(A^i \, x(t) - b^i \right), \label{eqn:g_i}
      \end{align}
      and $d$ vectors $R^i_1(t), \ldots, \, ~R^i_d(t)$, defined for all $j$ as
      \begin{align}\label{eqn:Rij}
        R^i_j(t)  & = \left(A^i\right)^T A^i k_j(t)  -  \left(\dfrac{1}{m}\right) e_j,
      \end{align}
      where $I$ denotes the $d$-dimensional identity matrix, and $e_j$ denotes the $j$-th column of the identity matrix $I$.
      \State The server updates the {\em pre-conditioner matrix} $K(t)$ to $K(t + 1)$ whose $j$-th column $k_j(t+1)$ is defined as
      \begin{align}
        k_j(t + 1) = k_j(t) - \alpha \textstyle \sum_{i=1}^m R^i_j(t). \label{eqn:kcol_update}
      \end{align}
      \State The server updates {\em the current estimate} $x(t)$ to  
      \begin{align}
        x(t+1) = x(t) - \delta K(t + 1) \textstyle \sum_{i=1}^m g^i(t). \label{eqn:x_update}
      \end{align}
    \EndFor
  \end{algorithmic}
\end{algorithm}

In this section, we present our key results on the robustness of the IPG method~\cite{chak2020ipc}, described in Algorithm~\ref{algo_1} for the noise-free case, in the presence of additive system noises; the {\em observation noise} and the {\em process noise}.\\

We introduce below some notation, our main assumption, and review a pivotal prior result. \\

\subsection{Notation, Assumption and Prior Results}
\label{sub:obs}
We define the {\em collective input matrix} $A$ and the {\em collective output vector} $b$ respectively as
\begin{align}
    A = \begin{bmatrix} (A^1)^T, \ldots, \, (A^m)^T \end{bmatrix}^T, ~ b  = \begin{bmatrix} (b^1)^T, \ldots, \, (b^m)^T \end{bmatrix}^T. \label{eqn:data_matrix}
\end{align}
Note that the matrix $\A$ is symmetric positive semi-definite and thus, has real non-negative eigenvalues. We let $\lambda_1 \geq \ldots \geq \lambda_d \geq 0$ denote the eigenvalues of $\A$ in order.

\begin{assumption}\label{assump:rank}
Assume that the matrix $\A$ is full rank.
\end{assumption}

Note that Assumption~\ref{assump:rank} holds true if and only if the matrix $\A$ is positive definite with $\lambda_d > 0$. Under Assumption 1, $\A$ is invertible, and we let $K^* = \left(\A\right)^{-1}$.
Note, from~\eqref{eqn:dist_opt}, that the Hessian of the aggregate cost function $\sum_{i = 1}^m F^i(x)$ is equal to $\A$ for all $x \in \R^d$. Thus, under Assumption 1 when $\A$ is positive definite, the aggregate cost function has a unique minimum point. Equivalently, the solution of the least squares problem defined by~\eqref{eqn:opt_1} is unique.\\




\input{meas_noise}

\input{sys_noise}

%% file: meas_noise.tex
We review below in Lemma~\ref{pro:pro_k} a prior result~\cite[Lemma~1]{chak2020ipc} that is pivotal for our key results presented later in this section. 
We let 
\begin{align}
    \varrho = (\lambda_1 - \lambda_d)/(\lambda_1 + \lambda_d). \label{eqn:varrho}
\end{align}
For each iteration $t$, recall the {\em pre-conditioner} matrix $K(t)$ in Algorithm~\ref{algo_1}, and define
\begin{align}
    \widetilde{K}(t) = K(t) - K^*, \quad \forall t \geq 0. \label{eqn:tilde_k}
\end{align}
Let $\widetilde{k}_j(t)$ denote the $j$-th column of $\widetilde{K}(t)$ and let $k_{j}^*$ denote the $j$-th column of the matrix $K^*$.
Lemma~\ref{pro:pro_k} below states sufficient conditions under which the sequence of the pre-conditioner matrices $\{K(t), ~ t = 0, \, 1, \ldots \}$, in Algorithm~\ref{algo_1}, converges linearly to $K^*$.

\begin{lemma}[Lemma~1 in~\cite{chak2020ipc}] \label{pro:pro_k}
Consider Algorithm~\ref{algo_1} with $\alpha \in (0,\frac{2}{\lambda_1})$. If Assumption 1 holds true then there exists a real value $\rho \in [\varrho, \, 1)$ such that, for all $j \in \{1, \ldots, \, d\}$,
\begin{align}
    \norm{\widetilde{k}_j(t+1)} \leq \rho \norm{\widetilde{k}_j(t)}, \quad \forall t \geq 0. \label{eqn:norm_tilde_K}
\end{align}
\end{lemma}
~


We first present in Section~\ref{sec:meas_noise} below the robustness of the IPG method against observation noise, followed by the robustness against process noise in Section~\ref{sec:sys_noise}.

\subsection{Robustness against Observation Noise}
\label{sec:meas_noise}

Based upon the literature~\cite{zhu2004class}, we model observation noise as follows. Each agent $i$ observes a corrupted output vector $b^{io}$, instead of the true output vector $b^{i}$. 
Specifically,
\begin{align}
    b^{io} = b^{i} + w_b^i, ~ i \in \{1,\ldots,m\}, \label{eqn:noisy_b_i}
\end{align}
where
$w_b^i \in \R^{n^i}$ is a random vector.
Let $\fixedexp{\cdot}$ denote the expectation of a function of the random vectors $\{w_b^i, ~ i=1,\ldots,m\}$. Let $\norm{\cdot}_1$ denote the $l^1$-norm~\cite{horn2012matrix}.
\begin{assumption}\label{assump:obsv_noise}
Assume that there exists $\eta < \infty$ such that 
\begin{align}
    \E \left[\norm{w_b^i}_1\right] \leq \eta, ~ i \in \{1,\ldots,m\}. \label{eqn:noise_bd_1}
\end{align}
\end{assumption}

In the presence of the above observation noise, in each iteration $t$ of Algorithm~\ref{algo_1}, each agent $i$ sends to the server a corrupted gradient $g^i(t)$ defined by~\eqref{eqn:gi_noisy} below, instead of~\eqref{eqn:g_i}. Specifically, for all $i$ and $t$,
\begin{align}
    g^i(t) = \left( A^i \right)^T \, \left(A^i \, x(t) - b^{io}\right). \label{eqn:gi_noisy}
\end{align}
Due to the above corruption in gradient computation, Algorithm~\ref{algo_1} no longer converges to an exact solution, defined by~\eqref{eqn:opt_1}, but rather to an approximation.\\

Theorem~\ref{thm:meas_noise} below presents a key result on the robustness of Algorithm~\ref{algo_1} against the above additive observation noise. Recall from Section~\ref{sub:obs} that under Assumption~1 the solution, denoted by $x^*$, of the regression problem~\eqref{eqn:opt_1} is unique. For each iteration $t$, we define the estimation error
\begin{align}
    z(t) = x(t) - x^*. \label{eqn:err}
\end{align}
For a matrix $M$, with columns $M_1, \ldots, M_d$, its Frobenius norm is defined to be $\norm{M}_F = \sqrt{\sum_{j=1}^d \norm{M_j}^2}$ \cite{horn2012matrix}. Recall the definition of $\varrho$ from~\eqref{eqn:varrho}.

\begin{theorem}\label{thm:meas_noise}
Consider Algorithm~\ref{algo_1} with parameters $\alpha < \frac{2}{\lambda_1}$ and $\delta \leq 1$, in the presence of additive observation noise defined in~\eqref{eqn:noisy_b_i} and the gradients for each iteration $t$, $\{g^i(t), i = 1, \ldots, \, m\}$, defined by~\eqref{eqn:gi_noisy}. If Assumptions~\ref{assump:rank}-\ref{assump:obsv_noise} hold then there exists $\rho \in [\varrho, \, 1)$ such that, for all $t \geq 0$, 
\begin{align}
    & \fixedexp{\norm{z(t+1)}} \leq 
    \left(1-\delta + \delta \lambda_1 \norm{\widetilde{K}(0)}_F \rho^{t+1} \right)\norm{z(t)}
    + \delta \eta m \sqrt{\lambda_1} \norm{\widetilde{K}(0)}_F \rho^{t+1} + \delta \eta m \sqrt{1/\lambda_d}. \label{eqn:rate_meas}
\end{align}
Additionally, $\lim_{t \rightarrow \infty} \fixedexp{\norm{z(t)}} \leq \delta \eta m \, \sqrt{1/\lambda_d}$.
\end{theorem}

Proof of Theorem~\ref{thm:meas_noise} is deferred to Section~\ref{pf:meas_noise}.\\

Since $\rho \in [0,1)$ and $\delta \in (0,1]$, Theorem~\ref{thm:meas_noise} implies that the IPG method under the influence of additive observation noise~\eqref{eqn:noisy_b_i} converges in expectation within a distance of $\delta \eta m \sqrt{1/\lambda_d}$ from the true solution $x^*$ of the regression problem defined in~\eqref{eqn:opt_1}.
Since the gradient-descent (GD) method (with constant step-sizes) is a special case of the IPG method with $K(t) = I$ for all $t$, from the proof of Theorem~\ref{thm:meas_noise} we obtain that the final error of GD
is bounded by $\delta \eta m \sqrt{\lambda_1}$.



%% file: sys_noise.tex
\subsection{Robustness against Process Noise}
\label{sec:sys_noise}
In the presence of process noise~\cite{gold1966effects, holi1993finite}, in each iteration $t$ of Algorithm~\ref{algo_1} (the IPG method in the noise-free case), the server computes corrupted values for both the pre-conditioner matrix $K(t)$ and the current estimate $x(t)$, described formally as follows. Specifically, in each iteration $t$, and for each $j \in \{1,  \dots, \, d \}$, instead of computing $k_j(t)$ (the $j$-th column of $K(t)$) accurately, the server computes
\begin{align}
    k^o_j(t) = k_j(t) + w^k_j(t), \label{eqn:noisy_k}
\end{align}
where $w^k_j(t) \in \R^d$. Similarly, instead of computing $x(t)$ accurately, the server computes a corrupted estimate 
\begin{align}
    x^o(t) = x(t) + w^x(t), \label{eqn:noisy_x}
\end{align}
where $w^x(t) \in \R^d$. Together, the vectors $\zeta(t) = \{w^x(t), w^k_j(t), j=1,\ldots,d\}$ are referred as additive process noise. 
For each iteration $t$, let $\serverexp{\cdot}$ denote the expectation of a function of the random vectors $\zeta(t)$.
Let $\expect{\cdot}$ denote the total expectation of a function of the random vectors $\{\zeta(0),\ldots,\zeta(t)\}$.
\begin{assumption}\label{assump:prc_noise}
Assume that the random vectors $\{w^x(t), w^k_j(t), j=1,\ldots,d\}$ are mutually independent for all $t$, and there exists $\omega < \infty$ such that for all $t, \,j$,
\begin{align}
    \serverexp{\norm{w^k_j(t)}_1} \leq \omega, \text{ and } \serverexp{\norm{w^x(t)}_1} \leq \omega. \label{eqn:noise_bd}
\end{align}
\end{assumption}

In the presence of above process noise, Algorithm~\ref{algo_1} is modified as follows. In each iteration $t$, each agent $i$ sends to the server a gradient $g^i(t)$ and $d$ vectors $R^i_1(t), \ldots, \, R^i_d(t)$ defined by~\eqref{eqn:gi_prc} and~\eqref{eqn:rij_prc} below, respectively, instead of~\eqref{eqn:g_i} and~\eqref{eqn:Rij}. Specifically, for all $i$ and $t$, 
\begin{align}
    g^i(t) = \left( A^i \right)^T \, \left(A^i \, x^o(t) - b^i\right),\label{eqn:gi_prc}
\end{align}
and, for $j = 1, \ldots, \, d$,
\begin{align}
    R^i_j(t) & = \left(A^i\right)^T A^i k_j^o(t)  -  \left(\dfrac{1}{m}\right) e_j. \label{eqn:rij_prc}
\end{align}
Similarly, in each iteration $t$, the server now computes $K(t+1)$ whose $j$-th column is defined as follows, instead of~\eqref{eqn:kcol_update}, for all $j \in \{1,\ldots, \,d\}$, 
\begin{align}
    k_j(t + 1) = k_j^o(t) - \alpha \textstyle \sum_{i=1}^m R^i_j(t). \label{eqn:kj_prc}
\end{align}
Recall that $k^o_j(t+1)$, defined by~\eqref{eqn:noisy_k}, is the corrupted value of $k_j(t+1)$.
Instead of~\eqref{eqn:x_update}, the updated estimate $x(t + 1)$ is defined by
\begin{align}
    x(t+1) = x^o(t) - \delta K^o(t + 1) \textstyle \sum_{i=1}^m g^i(t), ~ \forall t. \label{eqn:x_prc}
\end{align}
Recall that $x^o(t+1)$, defined by~\eqref{eqn:noisy_k}, is the corrupted value of $x(t+1)$.\\

To present our key result on the robustness of the IPG method, in Theorem~\ref{thm:sys_noise} below, against the above process noise, we introduce some notation. From Lemma~\ref{pro:pro_k}, recall that $\rho \in [\varrho,1)$. For each iteration $t$, we define
\begin{align}
    \widetilde{K}^o(t) & = \left[\widetilde{k}^o_1(t), \ldots, \, \widetilde{k}^o_d(t) \right] = K^o(t) - K^*, \text{ and } \label{eqn:obsrv_K}\\
    u(t) & = 1-\delta + \delta \lambda_1  \left(\rho^{t} \norm{\widetilde{K}(0)} + \omega\sqrt{d} \, \textstyle\sum_{i=0}^t \rho^i \right). \label{eqn:st}
\end{align}
Additionally, we let
\begin{align}
    \rho_{bd} = \frac{\norm{\widetilde{K}(0)}}{ \norm{\widetilde{K}(0)}+\omega\sqrt{d}} \, , ~ \text{ and } ~  \omega_{bd} = \frac{\left(1-\rho \right)}{
    \lambda_1 \sqrt{d}}. \label{eqn:rho_j}
\end{align}
~

Theorem~\ref{thm:sys_noise} below characterizes the convergence of Algorithm~\ref{algo_1} with modifications~\eqref{eqn:gi_prc}-\eqref{eqn:x_prc} in presence of process noise. Recall from definition~\eqref{eqn:err} that $z(t) = x(t) - x^*$ for all $t$.
Upon substituting $x(t)$ from~\eqref{eqn:noisy_x}, we obtain that
$z(t) = x^o(t) - w^x(t) -x^*$.
We let $z^0(t) = x^0(t) - x^*$. Thus, for each iteration $t$,
\begin{align}
    z^o(t) = z(t) + w^x(t). \label{eqn:noisy_z}
\end{align}

\begin{theorem} \label{thm:sys_noise}
Consider Algorithm~\ref{algo_1} with parameter $\alpha < \frac{2}{\lambda_1}$ and $\delta \leq 1$, in the presence of additive process noise defined in~\eqref{eqn:noisy_k}-\eqref{eqn:noisy_x}, and modifications defined in~\eqref{eqn:gi_prc}-\eqref{eqn:x_prc}.
If Assumptions~\ref{assump:rank} and~\ref{assump:prc_noise} hold then there exists $\rho \in [\varrho, \, 1)$ such that
\begin{align}
    \expect{\norm{z^o(t)}} & \leq \Pi_{k=1}^{t} u(k)\norm{z(0)}
    + \left(1+ u(t) + \ldots + \Pi_{k=1}^{t} u(k) \right)\omega, ~ \forall t. \label{eqn:x2}
\end{align}
Moreover, if
\begin{align}
    \rho & < \rho_{bd}, \text{ and }
    \omega < \omega_{bd}, \label{eqn:assump_2b}
\end{align}
then
$\lim_{t \rightarrow \infty} \expect{\norm{z^o(t)}} < \frac{\omega}{\delta \left(1 - \left(\omega/\omega_{bd}\right)\right)}$.
\end{theorem}

A proof of Theorem~\ref{thm:meas_noise} is deferred to Section~\ref{pf:sys_noise}.\\

Note that, from the proof of Theorem~\ref{thm:sys_noise} when $K(t)=I, \, \forall t$, the asymptotic estimation error of the traditional GD method under the above process noise is bounded by $\frac{\omega}{1-\norm{I-\delta\A}}$.\\

We present below formal proofs of our key results in Theorems~\ref{thm:meas_noise} and~\ref{thm:sys_noise} above -- however, {\bf the reader may proceed to Section \ref{sec:exp} without loss of continuity}.


%% file: proofs.tex
\section{PROOFS OF THEOREMS~\ref{thm:meas_noise} AND~\ref{thm:sys_noise}}
\label{sec:proofs}

We first present below a proof of Theorem~\ref{thm:meas_noise}. Then, in Section~\ref{pf:sys_noise} we prove Theorem~\ref{thm:sys_noise}.

\subsection{Proof of Theorem~\ref{thm:meas_noise}}
\label{pf:meas_noise}

From~\eqref{eqn:noisy_b_i}, recall the definition of corrupted output vector $b^{io}$, for each agent $i \in \{1,\ldots,m\}$, which are corrupted by the observation noise vector $w_b^i$. We concatenate the corrupted output vectors and the observation noise vector for all the agents to respectively define
\begin{align}
    b^o & = \begin{bmatrix} (b^{1o})^T, \ldots, \, (b^{mo})^T \end{bmatrix}^T, \label{eqn:bo_def} \\
    w_b & = \begin{bmatrix} (w_b^1)^T, \ldots, \, (w_b^m)^T \end{bmatrix}^T. \label{eqn:wb_def}
\end{align}
Upon substituting from~\eqref{eqn:noisy_b_i} in~\eqref{eqn:bo_def} we have
\begin{align*}
    b^o = \begin{bmatrix} (b^1)^T, \ldots, \, (b^m)^T \end{bmatrix}^T + \begin{bmatrix} (w_b^1)^T, \ldots, \, (w_b^m)^T \end{bmatrix}^T.
\end{align*}
Upon substituting from the definitions~\eqref{eqn:data_matrix} and~\eqref{eqn:wb_def} above we have
\begin{align}
    b^o = b + w_b. \label{eqn:noisy_b}
\end{align}
From the definition of vector norm we have~\cite{horn2012matrix}
\begin{align*}
    \norm{w_b} \leq \norm{w_b}_1 = \sum_{i=1}^m \norm{w_b^i}_1.
\end{align*}
From the definition of expectation, then we have
\begin{align*}
    \fixedexp{\norm{w_b}} \leq \sum_{i=1}^m \fixedexp{\norm{w_b^i}_1}.
\end{align*}
Under Assumption~\ref{assump:obsv_noise}, upon substituting above from~\eqref{eqn:noise_bd_1} we obtain that 
\begin{align}
    \fixedexp{\norm{w_b}} \leq \sum_{i=1}^m \fixedexp{\norm{w_b^i}_1} \leq \eta m. \label{eqn:noise_bd_2}
\end{align}

Consider an arbitrary iteration $t \geq 0$.
Upon substituting from~\eqref{eqn:gi_noisy} in~\eqref{eqn:x_update} we obtain that
\begin{align*}
    x(t+1) & = x(t) - \delta K(t+1) \left(\sum_{i=1}^m (A^i)^T A^i \right) x(t) - \left( \sum_{i=1}^m (A^i)^T b^{io} \right).
\end{align*}
As $A^T A = \sum_{i=1}^m (A^i)^T A^i$ and $A^T b^o = \sum_{i=1}^m (A^i)^T b^{io}$, the above implies that
\begin{align*}
    x(t+1) & = x(t) - \delta K(t+1) A^T (Ax(t)-b^o).
\end{align*}
Upon substituting above from~\eqref{eqn:noisy_b} we have
\begin{align}
    x(t+1) & = x(t) - \delta K(t+1) A^T (Ax(t) - b - w_b). \label{eqn:x_central}
\end{align}
Now, consider a point $x^* \in \arg \min_x \sum_{i = 1}^m F^i(x)$ defined by~\eqref{eqn:opt_1}. As $\nabla \sum_{i = 1}^m F^i(x) = A^T (A x - b)$ and the true value of the total gradient must vanish at the minimum point $x^*$ of~\eqref{eqn:opt_1}, we have
\begin{align}
    A^T (A x^* - b) = 0_d. \label{eqn:new_X*}
\end{align}
where $0_d$ denotes the origin of $\R^d$.
Recall from~\eqref{eqn:err} that $z(t) = x(t) - x^*$ for all $t$. Upon subtracting $x^*$ from both sides of~\eqref{eqn:x_central}, from the definition of $z(t)$ we have
\begin{align*}
    z(t+1) & = z(t) - \delta K(t+1) A^T (Ax(t) - b - w_b).
\end{align*}
Upon substituting above from~\eqref{eqn:new_X*}, we have
\begin{align}
    & z(t+1) = \left(I- \delta K(t+1)\A\right) z(t) + \delta K(t+1)A^T w_b. \label{eqn:z_mltp}  
\end{align}
Upon substituting above from~\eqref{eqn:tilde_k}, we have
\begin{align*}
    & z(t+1) = \left(I - \delta K^*\A\right) z(t) - \delta \widetilde{K}(t+1)\A z(t)
    + \delta K^*A^T w_b + \delta \widetilde{K}(t+1)A^T w_b.
\end{align*}
Upon using triangle inequality above and the definition $K^*\A = I$ (ref. Section~\ref{sub:obs}),
we obtain that
\begin{align}
    \norm{z(t+1)} & \leq \left(1-\delta \right) \norm{z(t)} + \delta \norm{\widetilde{K}(t+1)\A z(t)}
    + \delta \norm{K^*A^T w_b} + \delta \norm{\widetilde{K}(t+1)A^T w_b}.\label{eqn:z_mltp_2}
\end{align}
From the definition of induced 2-norm of a matrix~\cite{horn2012matrix} and expectation,~\eqref{eqn:z_mltp_2} implies that
\begin{align*}
    \fixedexp{\norm{z(t+1)}} \leq & \left(1-\delta + \delta \norm{\widetilde{K}(t+1)}\norm{\A}\right)\norm{z(t)} \\
    & +  \delta \norm{K^*A^T}\fixedexp{\norm{w_b}}
    + \delta \norm{\widetilde{K}(t+1)}\norm{A^T}\fixedexp{\norm{w_b}}.
\end{align*}
Upon substituting above from~\eqref{eqn:noise_bd_2} we obtain that
\begin{align}
    \fixedexp{\norm{z(t+1)}} \leq & \left(1-\delta + \delta \norm{\widetilde{K}(t+1)}\norm{\A}\right)\norm{z(t)} \nonumber\\
    & + \delta \norm{K^*A^T}\eta m 
    + \delta \norm{\widetilde{K}(t+1)}\norm{A^T}\eta m.\label{eqn:z_mltp_4}
\end{align}
Upon iterating the result in Lemma~\ref{pro:pro_k} from iteration $t$ to $0$, we have for all $j \in \{1, \ldots, \, d\}$,
\begin{align}
    \norm{\widetilde{k}_j(t)}^2 \leq \rho^{2t} \norm{\widetilde{k}_j(0)}^2. \label{eqn:k_norm_recur}
\end{align}
For a matrix $M$, with columns $M_1, \ldots, M_m$, its Frobenius norm is defined to be $\norm{M}_F = \sqrt{\sum_{j=1}^d \norm{M_j}^2}$. Note that $\norm{M} \leq \norm{M}_{F}$~\cite{horn2012matrix}.
Thus,
\begin{align*}
    \norm{\widetilde{K}(t)}^2 \leq \norm{\widetilde{K}(t)}_F^2 = \sum_{j=1}^d \norm{\widetilde{k}_j(t)}^2
\end{align*}
Upon substituting above from~\eqref{eqn:k_norm_recur} we have
\begin{align*}
    \norm{\widetilde{K}(t)}^2 \leq \rho^{2t} \sum_{j=1}^d \norm{\widetilde{k}_j(0)}^2 = \rho^{2t} \norm{\widetilde{K}(0)}_F^2,
\end{align*}
which means that
\begin{align}
    \norm{\widetilde{K}(t)} \leq \rho^{t} \norm{\widetilde{K}(0)}_F.
\end{align}
From basic Linear Algebra~\cite{horn2012matrix}, we know that $\norm{\A} = \lambda_1$ and $\norm{A^T} = \sqrt{\lambda_1}$.
Upon substituting these in~\eqref{eqn:z_mltp_4} we have
\begin{align}
    \fixedexp{\norm{z(t+1)}}  \leq & \delta \rho^{t+1} \norm{\widetilde{K}(0)}_F \left(\lambda_1\norm{z(t)} + \eta m\sqrt{ \lambda_1}\right) 
     + \left(1-\delta \right)\norm{z(t)} + \delta \norm{K^*A^T}\eta m. \label{eqn:zt_rhs_2}
\end{align}
From Singular Value Decomposition~\cite{horn2012matrix}, 
$$A = U S V^T,$$
where 
$S^T = \left[Diag\left(\sqrt{\lambda_1}, \ldots, \, \sqrt{\lambda_d} \right), ~ {\bf 0}_{d\times(\sum_{i = 1}^m n_i - d)}\right]$,
and the matrices $U, \, V$, respectively, constitutes of left and right orthonormal singular vectors of $A$. From above, $$(\A)^{-1}A^T  = V (S^T S)^{-1}S^T U^T.$$
Thus~\cite{horn2012matrix}, $$\norm{(\A)^{-1}A^T} = 1/\sqrt{\lambda_d}.$$
As $K^* = \left(\A\right)^{-1}$ (see Section~\ref{sub:obs}), the above implies that $$\norm{K^*A^T} = 1/\sqrt{\lambda_d}.$$ As discussed in Section~\ref{sub:obs}, Assumption~\ref{assump:rank} is equivalent to $\lambda_d > 0$. Upon substituting from above in~\eqref{eqn:zt_rhs_2} we obtain that
\begin{align*}
    \fixedexp{\norm{z(t+1)}} \leq 
    \left(1-\delta + \delta \lambda_1 \norm{\widetilde{K}(0)}_F \rho^{t+1} \right)\norm{z(t)}
    + \delta \eta m \sqrt{\lambda_1} \norm{\widetilde{K}(0)}_F \rho^{t+1} + \delta \eta m \sqrt{1/\lambda_d}.
\end{align*}
As $t$ is an arbitrary iteration, the above proves~\eqref{eqn:rate_meas}. \\

As $\rho \in [0,1)$ and $\delta \in (0, 1]$, there exists $T < \infty$ such that $\left(1-\delta + \delta \lambda_1 \norm{\widetilde{K}(0)}_F \rho^{t+1} \right) < 1$ for all $t\geq T$.
Thus, upon retracing~\eqref{eqn:rate_meas} from $t$ to $0$,
we have $\lim_{t \rightarrow \infty} \fixedexp{\norm{z(t)}} \leq \delta \eta m \, \sqrt{1/\lambda_d}$.



\subsection{Proof of Theorem~\ref{thm:sys_noise}}
\label{pf:sys_noise}
Similar to~\eqref{eqn:z_mltp} in Section~\ref{pf:meas_noise} above, for Algorithm~\ref{algo_1} with modifications~\eqref{eqn:gi_prc}-\eqref{eqn:x_prc}, we obtain that
\begin{align}
    z(t+1)  = \left(I- \delta K^o(t+1)\A\right) z^o(t), ~ \forall t. \label{eqn:above_1}
\end{align}
Upon substituting from~\eqref{eqn:obsrv_K} and~\eqref{eqn:noisy_z} in~\eqref{eqn:above_1} we have
\begin{align*}
    z(t+1) = \left(I-\delta K^* \A - \delta \widetilde{K}^o(t+1) \A\right) \left( z(t)+w^x(t) \right).
\end{align*}
Since $K^* \A = I$ (ref. Section~\ref{sub:obs}), from above we have
\begin{align*}
    z(t+1) = \left(\left(1-\delta \right)I - \delta \widetilde{K}^o(t+1) \A\right) \left( z(t)+w^x(t) \right).
\end{align*}
Using triangle inequality and the definition of induced matrix $2$-norm~\cite{horn2012matrix} above, we obtain that
\begin{align*}
    \norm{z(t+1)} \leq & \left(1-\delta + \delta \norm{\widetilde{K}^o(t+1)} \norm{\A} \right) \left(\norm{z(t)}+\norm{w^x(t)} \right) \\
    = & \left(1-\delta + \delta \lambda_1 \norm{\widetilde{K}^o(t+1)} \right) \left(\norm{z(t)}+\norm{w^x(t)} \right).
\end{align*}
From Assumption~\ref{assump:prc_noise}, since the random variables $\{w^x(t), w^k_j(t), j=1,\ldots,d\}$ are mutually independent for all $t$, the above implies that
\begin{align}
    \expect{\norm{z(t+1)}} \leq & \left(1-\delta + \delta \lambda_1 \expect{\norm{\widetilde{K}^o(t+1)}}\right) \left(\norm{z(t)}+\expect{\norm{w^x(t)}} \right). \label{eqn:zt_full_1}
\end{align}

Define the matrix $$W^k(t) = \begin{bmatrix} w^k_1(t), \ldots, \, w^k_d(t) \end{bmatrix}.$$
From the definition of $K^*$ in~\eqref{eqn:obsrv_K} we get
\begin{align*}
    \widetilde{K}^o(t) = \begin{bmatrix} k^o_1(t), \ldots, \, k^o_d(t) \end{bmatrix} - K^*.
\end{align*}
Recall the definition of additive process noise in~\eqref{eqn:noisy_k}.
Upon substituting from~\eqref{eqn:noisy_k} in~\eqref{eqn:obsrv_K} we have 
\begin{align*}
    \widetilde{K}^o(t) & = \begin{bmatrix} k_1(t), \ldots, \, k_d(t) \end{bmatrix} + \begin{bmatrix} w^k_1(t), \ldots, \, w^k_d(t) \end{bmatrix} - K^* = K(t) + W^k(t) - K^* = \left(K(t) - K^* \right) + W^k(t).
\end{align*}
Upon substituting above from~\eqref{eqn:tilde_k} we obtain that 
$$\widetilde{K}^o(t) = \widetilde{K}(t) + W^k(t).$$
Using triangle inequality above we have
\begin{align*}
    \norm{\widetilde{K}^o(t)} \leq \norm{\widetilde{K}(t)} + \norm{W^k(t)},
\end{align*}
which implies that
\begin{align}
    \expect{\norm{\widetilde{K}^o(t)}} \leq \norm{\widetilde{K}(t)} + \expect{\norm{W^k(t)}}. \label{eqn:kt_split}
\end{align}
From the definition of matrix norms~\cite{horn2012matrix},
\begin{align*}
    \norm{W^k(t)} \leq \sqrt{d}~\norm{W^k(t)}_1 = \sqrt{d}~\max_{j=1,\ldots,d} \norm{w^k_j}_1.
\end{align*}
From the definition of expectation then we have
\begin{align*}
    \expect{\norm{W^k(t)}} \leq \sqrt{d}~\max_{j=1,\ldots,d}\expect{\norm{w^k_j}_1}.
\end{align*}
Under Assumption~\ref{assump:prc_noise}, upon substituting above from~\eqref{eqn:noise_bd} we obtain that
$$\expect{\norm{W^k(t)}} \leq  \omega\sqrt{d}.$$
Upon substituting from above in~\eqref{eqn:kt_split} we get
\begin{align}
    \expect{\norm{\widetilde{K}^o(t)}} \leq \norm{\widetilde{K}(t)} + \omega\sqrt{d}. \label{eqn:kt_deter}
\end{align}
Instead of each column $\widetilde{k}_j(t)$ in Lemma~\ref{pro:pro_k} if we consider $\widetilde{K}(t)$ then following the proof of Lemma~\ref{pro:pro_k} for $\widetilde{K}_j(t)$, we have a similar result as~\eqref{eqn:norm_tilde_K}: 
\begin{align*}
    \norm{\widetilde{K}(t+1)} \leq \rho \norm{\widetilde{K}(t)}, ~ \forall t\geq 0.
\end{align*}
Upon substituting from above in~\eqref{eqn:kt_deter} we get
\begin{align*}
    \expect{\norm{\widetilde{K}^o(t)}} \leq \rho \norm{\widetilde{K}(t-1)} + \omega\sqrt{d} \leq \rho^{t}\norm{\widetilde{K}(0)} + \omega\sqrt{d}\sum_{i=0}^{t} \rho^i.
\end{align*}
Upon substituting above from~\eqref{eqn:st} we obtain that
\begin{align*}
    \expect{\norm{\widetilde{K}^o(t)}} \leq \left(u(t)-1+\delta\right) / \left(\delta\lambda_1\right).
\end{align*}
Upon substituting from above in~\eqref{eqn:zt_full_1} we get
\begin{align*}
    \expect{\norm{z(t+1)}}
    & \leq u(t+1)\left(\norm{z(t)}+\expect{\norm{w^x(t)}} \right) \leq u(t+1) \left(\norm{z(t)}+\expect{\norm{w^x(t)}_1} \right).
\end{align*}
Under Assumption~\ref{assump:prc_noise}, upon substituting above from~\eqref{eqn:noise_bd} we obtain that
\begin{align}
    \expect{\norm{z(t+1)}}
    & \leq u(t+1) \left(\norm{z(t)}+\omega \right) \nonumber \\
    & \leq \Pi_{k=1}^{t+1} u(k) \norm{z(0)} + \left(u(t+1) + u(t+1) u(t) + \ldots + \Pi_{k=1}^{t+1} u(k) \right)\omega. \label{eqn:znot_3}
\end{align}
From~\eqref{eqn:noisy_z} and~\eqref{eqn:noise_bd} we have $$\expect{\norm{z^o(t+1)}} \leq \expect{\norm{z(t+1)}} + \omega.$$
Substituting from~\eqref{eqn:znot_3} in the R.H.S.~above proves~\eqref{eqn:x2}.\\

As $\rho < 1$, from the definition of $u(t)$ in~\eqref{eqn:st}, we obtain the limiting value of $u(t)$ as
\begin{align*}
    \lim_{t \rightarrow \infty} u(t) =  1-\delta + \delta \frac{ \lambda_1  \omega \sqrt{d}}{1-\rho}.
\end{align*}
Upon substituting above from the definition of $\omega_{bd}$ in~\eqref{eqn:rho_j} we get
\begin{align}
    \lim_{t \rightarrow \infty} u(t) =  1-\delta + \delta\frac{\omega}{\omega_{bd}}. \label{eqn:r_lim}
\end{align}
The condition $\omega < \omega_{bd}$ in~\eqref{eqn:assump_2b} then implies that
\begin{align*}
    \lim_{t \rightarrow \infty} u(t) < 1. 
\end{align*}
From~\eqref{eqn:st}, we further have that
\begin{align*}
    u(t) - u(t-1) = \delta \lambda_1 \rho^{t-1} \left(\rho \left(\omega\sqrt{d}+\norm{\widetilde{K}(0)}\right) - \norm{\widetilde{K}(0)} \right), ~ \forall t \geq 1.
\end{align*}
Upon substituting above from the definition of $\rho_{bd}$ in~\eqref{eqn:rho_j} we have
\begin{align*}
    & u(t) - u(t-1) =  \delta \lambda_1 \rho^{t-1} \norm{\widetilde{K}(0)} \left( (\rho/\rho_{bd}) - 1\right), ~ \forall t \geq 1.
\end{align*}
The above, in conjunction with the condition~\eqref{eqn:assump_2b}, implies that
$$u(t) < u(t-1), ~ \forall t \geq 1.$$
The limit in~\eqref{eqn:r_lim}, in conjunction with the fact that $u(t)$ is non-negative for all $t$,
implies that there exists $\tau < \infty$ such that $0 \leq u(t) < 1$ for all $t \geq \tau$.
Thus,
\begin{align}
    \lim_{t \rightarrow \infty} \Pi_{k=1}^t u(k) & = 0, \, \text{and} \label{eqn:x3} \\
    \lim_{t \rightarrow \infty} \left(1+ u(t) + \ldots + \Pi_{k=1}^{t} u(k) \right) & < \dfrac{1}{1 - u(\tau)}. \label{eqn:x4}
\end{align}
Substituting from~\eqref{eqn:x3} and~\eqref{eqn:x4} into~\eqref{eqn:x2}, we obtain that
\begin{align}
    \lim_{t \rightarrow \infty} \expect{\norm{z^o(t)}}  < \omega/(1 - u(\tau)). \label{eqn:sse}
\end{align}
Since $\{u(t)\}$ is a strictly decreasing sequence, if~\eqref{eqn:sse} holds true for some $\tau$ satisfying $u(\tau) < 1$ then it also holds true for $\tau+1, \, \ldots, \infty$. Thus,~\eqref{eqn:sse} implies that 
$$\lim_{t \rightarrow \infty} \expect{\norm{z^o(t)}}  < \omega/(1 - u(\infty)).$$ 
Upon substituting $u(\infty)$ above from~\eqref{eqn:r_lim} we obtain that
$$\lim_{t \rightarrow \infty} \expect{\norm{z^o(t)}} < \frac{\omega}{\delta \left(1 - \left(\omega/\omega_{bd}\right)\right)}.$$ Hence, the proof.


%% file: simulations.tex
\section{EXPERIMENTAL RESULTS}
\label{sec:exp}


\begin{table*}[htb!]
\caption{\it The parameters used in different algorithms for their minimum convergence rate~\cite{chakrabarti2020iterative,azizan2019distributed}.}
\begin{center}
\begin{tabular}{|p{1.2cm}||p{2.3cm}|p{1.5cm}|p{3cm}|p{3cm}|p{2.5cm}|}
\hline
Dataset & {\bf IPG}~\cite{chakrabarti2020iterative} & GD~\cite{azizan2019distributed} & NAG~\cite{azizan2019distributed} & HBM~\cite{azizan2019distributed} & APC~\cite{azizan2019distributed} \\
\hline
\hline
ash608 & $\alpha = 0.1163, \, \delta = 1$ & $\alpha = 0.1163$ & $\alpha = 0.08, \, \beta = 0.5$ & $\alpha = 0.15, \, \beta = 0.29$ & $\gamma = 1.02, \, \eta =5.27$ \\
\hline
gr\_30\_30 & $\alpha = 0.014, \, \delta = 1$ & $\alpha = 0.014$ & $\alpha = 0.009, \, \beta =0.99$ & $\alpha = 0.03, \, \beta =0.98$ & $\gamma = 1.09, \, \eta =12.8$ \\
\hline
\end{tabular}
\end{center}
\label{tab:parameters}
\end{table*}

\begin{figure*}[htb!]
\centering
\begin{subfigure}{.5\textwidth}
  \begin{center}
  \includegraphics[width = \textwidth]{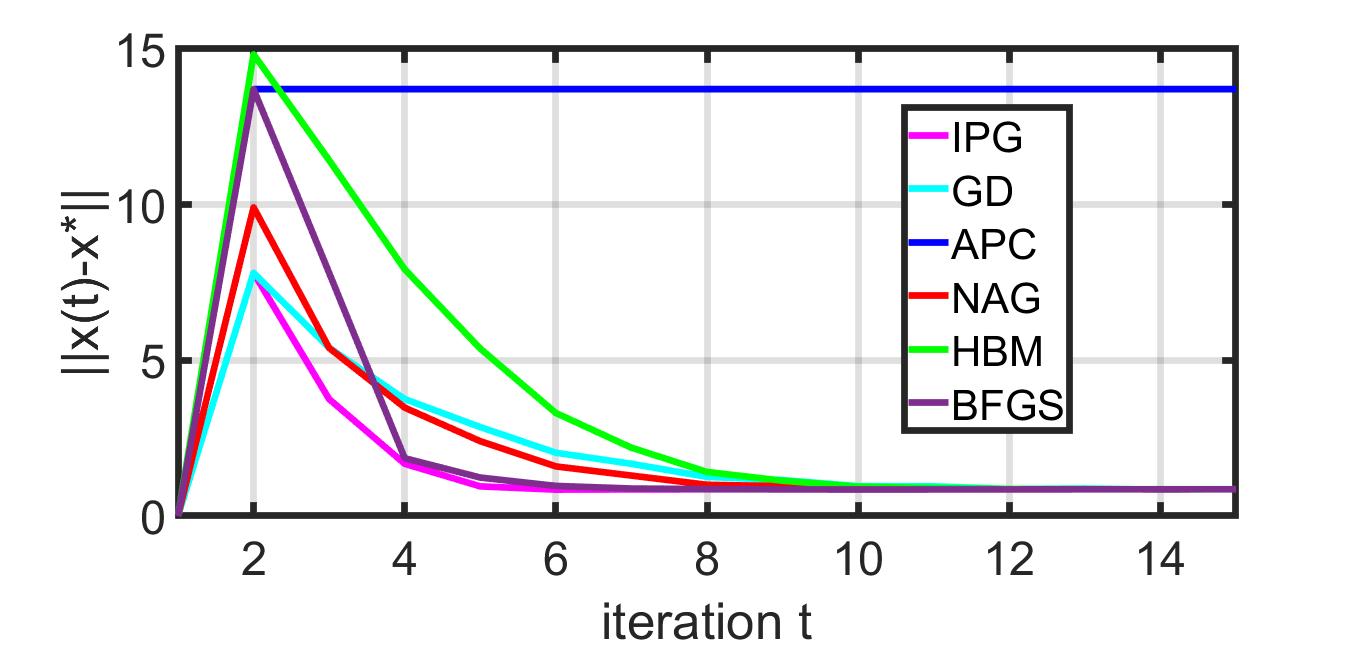}
  \caption{``ash608'': observation noise}
  \end{center}
\end{subfigure}%
\begin{subfigure}{.5\textwidth}
  \begin{center}
  \includegraphics[width = \textwidth]{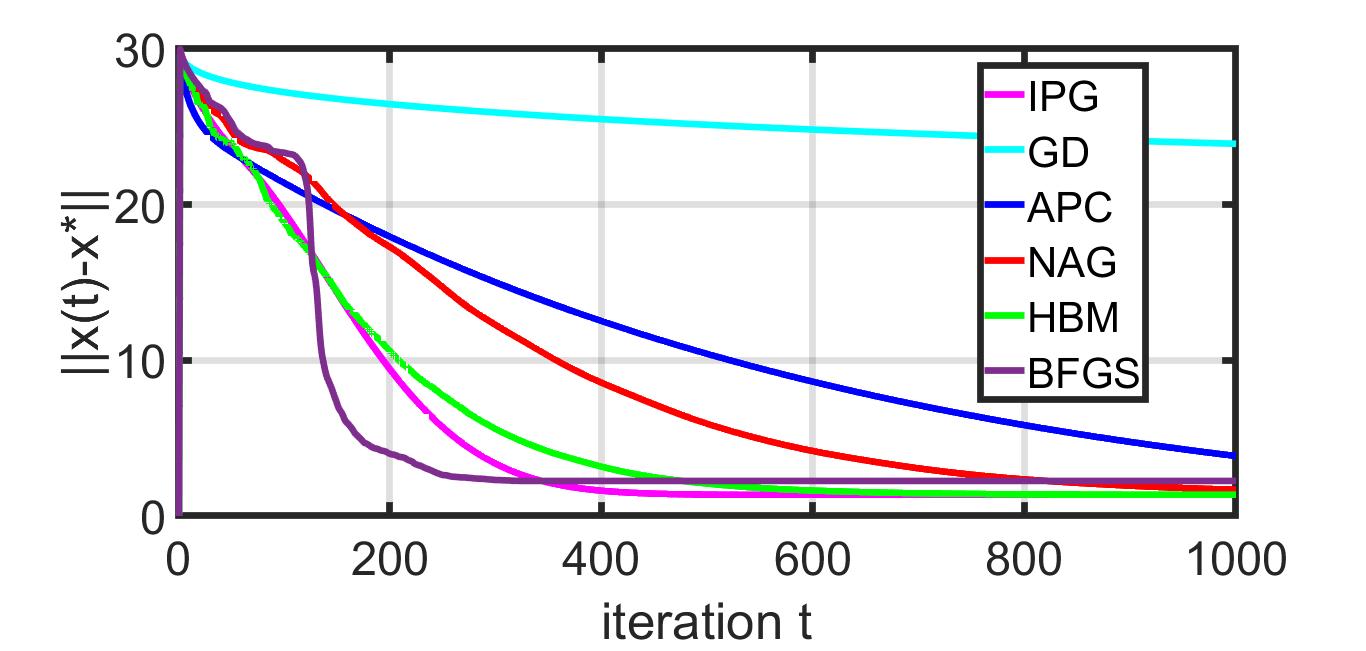}
  \caption{``gr\_30\_30'': observation noise}
  \end{center}
\end{subfigure}
\begin{subfigure}{.5\textwidth}
  \begin{center}
  \includegraphics[width = \textwidth]{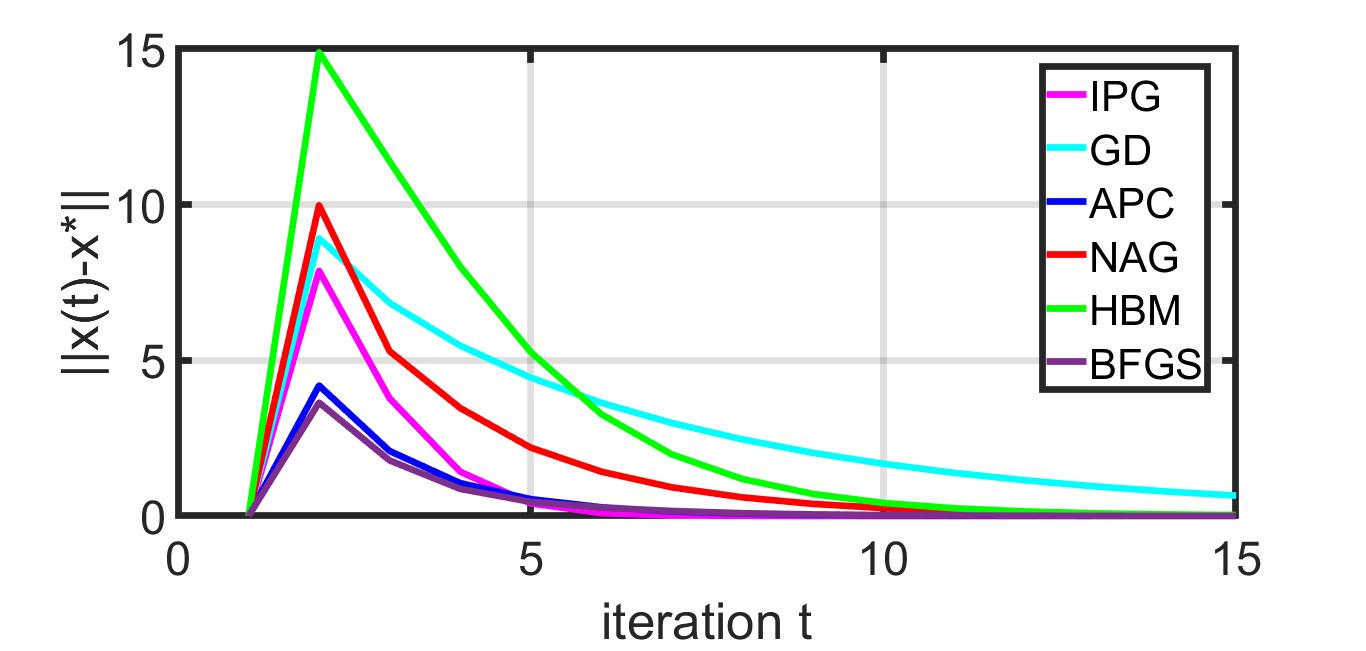}
  \caption{``ash608'': process noise}
  \end{center}
\end{subfigure}%
\begin{subfigure}{.5\textwidth}
  \begin{center}
  \includegraphics[width = \textwidth]{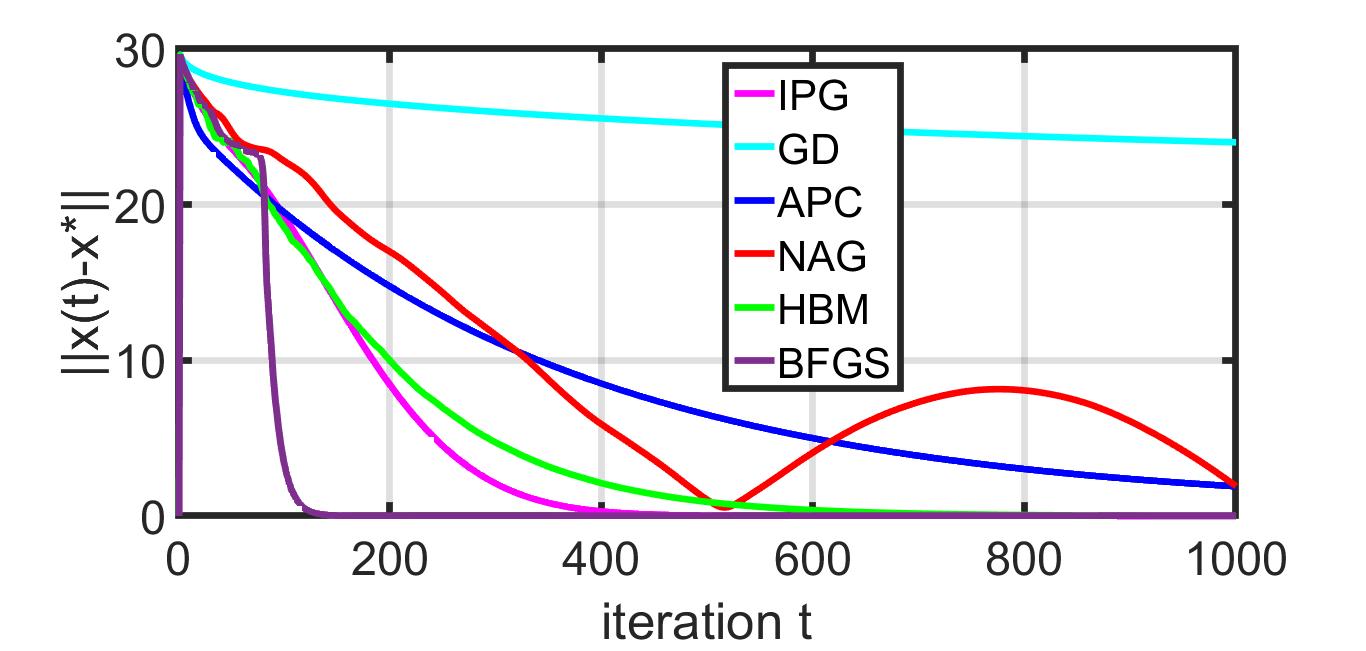}
  \caption{``gr\_30\_30'': process noise}
  \end{center}
\end{subfigure}
\caption{\footnotesize{\it Temporal evolution of estimation error $\norm{x(t)-x^*}$ in presence of additive observation noise (shown in the first row) and process noise (shown in the second row), for the different algorithms represented by different colors. For algorithms IPG, GD, NAG, HBM, and BFGS, $x(0) = [0,\ldots,0]^T$. Additionally, for IPG, $K(0) = O_{d \times d}$, and for BFGS, $M(0) = I$ the identity matrix. APC is initialized as per the instructions in~\cite{azizan2019distributed}. The other parameters are enlisted in Table~\ref{tab:parameters}.}}
\label{fig:meas_noise}
\end{figure*}

\begin{table*}[htb!]
\caption{\it Comparisons between the final estimation errors $\lim_{t \rightarrow \infty} \norm{x(t)-x^*}$ for different algorithms.}
\begin{center}
\begin{tabular}{|p{2cm}|p{1.2cm}|p{1.5cm}||p{0.6cm}|p{1.42cm}|p{1.42cm}|p{1.42cm}|p{1.42cm}|p{1.42cm}|}
\hline
Noise type & Dataset & Noise level & {\bf IPG} & GD & NAG & HBM & APC & BFGS \\ 
\hline
\hline
Observation & ash608 & $\eta=8.23$ & $0.86$ & $0.86$ & $0.86$ & $0.86$ & $13.71$ & $0.86$ \\
noise & gr\_30\_30 & $\eta =7.21$ & $1.35$ & $2.05$ & $1.35$ & $1.35$ & $1.82$ & $2.25$ \\
\hline
\multirow{ 2}{*}{Process noise} & ash608 & $\omega=9.3 \times 10^{-3}$ & $0$ & $3.46 \times 10^{-4}$ & $9.21 \times 10^{-4}$ & $10^{-4}$ & $4.9 \times 10^{-4}$ & $\infty$ \\
& gr\_30\_30 & $\omega=4.5 \times 10^{-2}$ &  $0$ & $7.68$ & $1.86$ & $8.5 \times 10^{-3}$ & $3.72$ & $1.49 \times 10^{-2}$ \\
\hline
\end{tabular}
\end{center}
\label{tab:sse}
\end{table*}


We consider $2$ benchmark input matrices $A$ (see~\eqref{eqn:data_matrix}), namely \textit{``ash608''} and \textit{``gr\_30\_30''}, from the SuiteSparse Matrix Collection (\underline{http://sparse.tamu.edu/}). The true value of the collective output is $b=Ax^*$ where $x^*$ is a $d$-dimensional vector with all entries equal to $1$. The rows of $(A,b)$ are distributed among $m=10$ agents. For both these datasets, Assumption 1 holds true.\\

{\bfseries Observation Noise}:
We add {\em uniformly} distributed random noise vectors from $(-0.25,0.25)$ and $(-0.15,0.15)$, respectively, to the true output vectors of datasets \textit{``ash608''} and \textit{``gr\_30\_30''}. The algorithm parameters are chosen such that each algorithm achieves its minimum convergence rate~\cite{chakrabarti2020iterative,azizan2019distributed} (ref. Table~\ref{tab:parameters}). Each iterative algorithm is run until the changes in its subsequent estimates 
is less than $10^{-4}$ over $20$ consecutive iterations. We note the final estimation errors of the different algorithms in Table~\ref{tab:sse}. We observe that the final estimation error of the IPG method is either comparable or favourable to all the other algorithms, for each dataset. \\

{\bfseries Process noise}:
We simulate the algorithms by adding noise to the iterated variables. For the algorithms GD, NAG, HBM, and IPG, the process noise has been generated by rounding-off each entries of all the iterated variables in the respective algorithms to {\em four} decimal places. However, the rounding-off does not generate same values of noise $w$ for the APC and BFGS algorithms. Therefore, for APC, we add {\em uniformly} distributed random numbers in the range $(0,5\times10^{-5})$ for both the datasets, and similarly, for BFGS, we add {\em uniformly} distributed random numbers in the range $(0,9 \times 10^{-5})$ and $(0,2 \times 10^{-6})$ respectively for the datasets \textit{``ash608''} and \textit{``gr\_30\_30''}.
The final estimation errors of different algorithms are noted in Table~\ref{tab:sse}. We observe that the final error for IPG is less than all the other algorithms. Also, we observe that the estimation error for the BFGS algorithm on dataset \textit{``ash608''} grows unbounded after $360$ iterations. The cause for this instability is the violation of the non-singularity of the approximated Hessian matrix, which is a necessary condition for the convergence of BFGS~\cite{kelley1999iterative}.

%% file: summary.tex
\section{SUMMARY}

In this paper, we solve the distributed noisy linear regression problem in server-based network architecture. We have considered practical settings with additive system noises: either observation noise or process noise. Our contribution has been analyzing the recently proposed Iteratively Pre-Conditioned Gradient-descent (IPG) algorithm's robustness against such independent system noises whose magnitudes ($l^1$-norm) are bounded in expectation. The experimental results have reinforced our claim on the IPG method's superior accuracy compared to other state-of-the-art distributed algorithms when subjected to system noises.